\documentclass[12pt]{article}
\usepackage[cp1251]{inputenc}
\usepackage{amsmath}
\usepackage{amssymb}
\usepackage{mathrsfs}
\usepackage{amsfonts, dsfont, srcltx}

\newtheorem{lemma}{Lemma}
\newtheorem{theorem}{Theorem}

\begin{document}

\begin{center}
{\Large  Invertibility In the Sense of Ehrenpreis
}
\end{center}

\begin{center}
N. F.~{Abuzyarova}
\end{center}

[E-mail: ]{abnatf@gmail.com} 

{Bashkir State
University,  Zaki Validi street 32, Ufa, 450076, Russia}

\setcounter{page}{1}

\begin{center}
In this paper we consider   zero sets of entire functions belonging to the Schwartz algebra.
This algebra is defined as the Fourier-Laplace transform 
image of the space of all distributions compactly supported on the real line.  
We study the conditions under which given complex sequence forms  zero set of some invertible in
the sense of Ehrenpreis element of the Schwartz algebra (slowly decreasing function).
\end{center}

{30D15; 30E05; 42A38; 46F05}


\section{Introduction}

We will use the  following common notations.
$\mathcal D =C_0^{\infty} (\mathbb R)$ denotes the space of test functions,
$\mathcal D '=(C_0^{\infty} (\mathbb R))'$. 
$\mathcal E=C^{\infty} (\mathbb R )$ is the Schwartz space equipped with its standard metrizable topology,
   $\mathcal E'=(C^{\infty} (\mathbb R))'$ is its strong dual space consisting of all distributions compactly supported on the real line.
    
		Given $S\in\mathcal E',$
		its Fourier-Laplace transform $\mathcal F (S)$ is defined by the formula
	$$\mathcal F (S)=S(e^{-\mathrm{i}tz}).$$ 
	
	It is well-known that the image $\mathcal P=\mathcal F (\mathcal E')$ is
	the linear space of all entire functions of exponential type having polynomial growth along the real axis \cite[Theorem 7.3.1]{Horm}.
	Being equipped with the topology induced from $\mathcal E'$, $\mathcal P$ becomes a locally convex space of the type $(LN^*)$; 
	and $\mathcal P$ is also the topological algebra ({\sl Schwartz algebra}).

Recall that  $S\in\mathcal E'$ is {\sl invertible} distribution \cite{Ehren}
if
$$
S*\mathcal E=\mathcal E, 
$$
$$
S*\mathcal D' =\mathcal D',
$$
where $*$ denotes the convolution.

$S$ is invertible if and only if the primary ideal algebraically generated by $\varphi =\mathcal F (S)$ in  $\mathcal P$
is closed. It is also equivalent to the  <<division theorem>>:
{\it  $\Phi\in\mathcal P$ and $\Phi/\varphi\in Hol (\mathbb C)$ imply  $\Phi/\varphi\in\mathcal P.$ }

We  call the function $\varphi\in\mathcal P$
{\sl invertible in the sense of Ehrenpreis} if  $S=\mathcal F^{-1} (\varphi)$ is invertible.

Due to the results of the paper \cite{Ehren} (Theorem I, Theorem 2.2, Proposition 2.7), 
$\varphi \in\mathcal P$ is invertible in the sense of Ehrenpreis if and only if this function is
{\sl slowly decreasing}. It means that
 there exists $ a>0$ with the property
\begin{equation}
\forall x\in\mathbb R\ \exists x'\in\mathbb R: \ |x-x'|\le a\mathrm{ln}\, (a+|x|), \
|\varphi (x')|\ge (a+|x'|)^{-a} .
\label{treb}
\end{equation}

Below, we formulate two other equivalent definitions of slow decrease property 
 for function $\varphi \in\mathcal P$ (\cite[Section  3]{Ehren}):

1) there exists $ a>0$ such that  
\begin{equation}
\forall x\in\mathbb R\ \exists z'\in\mathbb C: \ |x-z'|\le a\mathrm{ln}\, (a+|x|), \
|\varphi (z')|\ge (a+|z'|)^{-a} ;
\label{treb-2}
\end{equation}

2) there exists $ a>0$ such that  $\forall x\in\mathbb R$ we can find a circumference $C_x$ of 
radius not greater than $a\mathrm{ln}\, (a+|x|)$ with   $x$  laying inside $C_x$ 
and 
$$
|\varphi (z)|\ge (a+|z|)^{-a} \ \forall z\in C_x .
$$

As a rule, 
given a class of entire functions defined by some growth restrictions, it is interesting and useful
to study the properties of their zero sets of these functions.
In particular, it is related to
 investigating  the necessary conditions, the sufficient conditions and the criteria
of invertibility in the sense of Ehrenpreis for the function $\varphi\in\mathcal P$ 
in terms of some  characteristics of its zero set.

There are well-known connections of
invertible in the sense of Ehrenpreis functions with 
question of surjectivity for the convolution operator in $\mathcal E$ and
 <<fundamental principle>> for the set of solutions of homogeneous convolution equation (or system of such equations)
 (see survey \cite{Ber-St}).

For instance,
if $S\in\mathcal E'$ 
then the zero set of 
$\varphi=\mathcal F (S)$ up to the multiplier $(-\mathrm{i})$
equals the spectrum of the solution subspace   
   for  the homogeneous convolution equation
$S*f=0.$ 
This subspace admits spectral synthesis \cite{Horm-2}.
Moreover, if $\varphi$ is invertible in the sense of Ehrenpreis then
each solution $f$ is represented as a series of exponential monomials
converging with respect to the topology of $\mathcal E$ 
after some groupings of its members
\cite[Theorem 3.1]{Ehren}.

It is natural to expect the similar role of invertible in the sense of Ehrenpreis functions
in questions of representation  and <<fundamental principle>>
for differential-invariant subspaces of the space $\mathcal E(a;b):=C^{\infty} (a;b),$
where $(a;b)$ is finite or infinite interval of the real axis.
Spectral synthesis problem for the differentiation operator in $\mathcal E(a;b)$
has studied since 2008. Most of the known results are contained in 
  \cite{Al-Kor}--\cite{Bar-Bel}.

In this paper, we obtain some necessary conditions and  common criteria of invertibility in the sense of Ehrenpreis 
for functions $\varphi\in\mathcal P$ which zero sets
 $\{\lambda_j\}$ 
satisfy the relation
$
\mathrm{Im}\, \lambda_j = O(\mathrm{ln} \, |\lambda_j|)
$
as $j\to\infty.$

In Section 2, we prove Lemma \ref{lem-1}.
This lemma allows to study only   functions which zeros are real.
Then, we  prove two 
assertions on the necessary conditions of invertibility in the sense of Ehrenpreis
  (Lemma \ref{lem-2} and Theorem \ref{tm-1}).
 
Section 3 contains  the necessary and sufficient conditions under which
 given  even entire function of  exponential type with  real zeros 
is an invertible in the sense of Ehrenpreis element of the algebra $\mathcal P$  (Theorem \ref{tm-2}).
  We also prove the common criterion of invertibility in the sense of Ehrenpreis
for an arbitrary function $\varphi\in\mathcal P$ which zeros are real
(Theorem \ref{tm-3}).

\section{Necessary conditions of  invertibility in the sense of Ehrenpreis}

Let $\mathcal M=\{\mu_j\},$ $\mu_j=\alpha_j+\mathrm{i}\beta_j \subset\mathbb C,$
$ 
0<|\mu_1|\le |\mu_2|\le\dots,
$
be such that
 $\beta_j=O(\mathrm{ln}\,|\mu_j|)$ as $j\to\infty ,$
 and
\begin{equation}
\psi(z)=\lim_{R\to\infty} \prod_{|\mu_j|\le R}\left( 1-\frac{z}{\mu_j}\right)
\label{fi-lam}
\end{equation}
is an entire function of exponential type

\begin{lemma}
The function $\psi \in\mathcal P$ 
and it is invertible in the sense of Ehrenpreis 
if and only if  
the function
\begin{equation}
\psi_1(z)=\lim_{R\to\infty} \prod_{|\alpha_j|\le R}\left( 1-\frac{z}{\alpha_j}\right) 
\label{fi-alf}
\end{equation}
is  invertible in the sense of Ehrenpreis  element of the algebra $\mathcal P$.
\label{lem-1}
\end{lemma}

{\bf Proof.}
First, notice that either both sequences, $\{\mu_j\}$ and $\{\alpha_j\}$, have   densities or not;
and either both series, $\sum\mu_j^{-1}$ and $\sum \alpha_j^{-1}$, converge or not.
It follows that either both formulas, (\ref{fi-lam}) and (\ref{fi-alf}), define  entire functions of the same exponential type or not.

For each multiplier in the right-hand side of 
(\ref{fi-alf}), we have 
$$
\left| 1-\frac{x}{\alpha_j}\right|\le\left| 1-\frac{x}{\mu_j}\right|\left( 1+\frac{\beta_j^2}{\alpha_j^2}\right)^{1/2},
\quad x\in\mathbb R.
$$
Hence, if $\psi$ is an entire function then
\begin{equation}
\mathrm{ln}\, |\psi_1 (x)|\le \mathrm{ln}\, |\psi (x)| + O(1),  \quad x\in\mathbb R.
\label{fi-fi-1}
\end{equation}
By this estimate, we get 
\begin{equation}
\psi\in \mathcal P\Longrightarrow \psi_1\in \mathcal P,
\label{i-1}
\end{equation}
and if $\psi\in\mathcal P$ then
\begin{multline}
\psi_1\ \text{is invertible in the sense of Ehrenpreis in} \ \mathcal P \Longrightarrow 
\\ \Longrightarrow 
\psi\ \text{is invertible in the sense of Ehrenpreis in} \ \mathcal P.
\label{i-2}
\end{multline}

Further, set
$$
\mathcal M^+=\{\mu_j:\ \beta_j\ge 0\},\quad \mathcal M^{-} =\mathcal M\setminus\mathcal M^+ ,
$$
$$
\psi^+ (z)=\lim_{R\to\infty}  \prod_{\substack{|\mu_j|\le R\\ \mu_j\in\mathcal M^-}}
\left( 1-\frac{z}{\mu_j}\right)\cdot\prod_{\substack{|\mu_j|\le R\\ \mu_j\in\mathcal M^+}}
\left( 1-\frac{z}{\alpha_j}\right) .
$$

Without loss of the generality ({\sl wlog}), we may assume that $|\alpha_j|>1,$  $j\in\mathbb N.$
Let $M_0>0$ be such that
$$
|\beta_j|\le M_0\mathrm{ln}\, |\alpha_j|,\  j=1,2,\dots
$$
It is easy to see that
$$
\left| 1-\frac{z}{\mu_j}\right|\le\left| 1-\frac{z}{\alpha_j}\right|,
$$
where 
$z=x+2\mathrm{i} M_0\mathrm{ln}\, |x| ,$  $|x|>2$, 
 $\mu_j\in\mathcal M^+$ and $|\alpha_j|\le x^4$.

\noindent
It is also true that
$$
\left| 1-\frac{z}{\mu_j}\right|\le\left| 1-\frac{z}{\alpha_j}\right|\cdot
\left( 1+\frac{4M^2_0\mathrm{ln}^2\,|\alpha_j|}{|\alpha_j|^2}\right)^{1/2}
$$
where 
$z=x+2\mathrm{i} M_0\mathrm{ln}\, |x| ,$  $|x|>2$, 
 $\mu_j\in\mathcal M^+$ and $|\alpha_j|> x^4$.

From these inequalities  we derive the estimate 
 \begin{equation}
\mathrm{ln}\, |\psi (z)|\le \mathrm{ln}\, |\psi^+ (z)| +O(1),  \quad z=x+2\mathrm{i} M_0\mathrm{ln}\,  |x|
\ \text{as } \
|x|\to \infty .
\label{e-1}
\end{equation}

By the same argument, we obtain 
\begin{equation}
\mathrm{ln}\, |\psi^+ (z)|\le \mathrm{ln}\, |\psi_1 (z)| +O(1), 
\quad z=x-2\mathrm{i} M_0\mathrm{ln}\,  |x|
\ \text{as } \
|x|\to \infty .
\label{e-2}
\end{equation}

By
(\ref{i-1})--(\ref{e-2}), 
taking into account that the functions $\psi,$  $\psi_1,$ $\psi^+$ have the  same exponential type
and applying the Phragmen-Lindel\"of principle,
we get
\begin{equation}
\psi_1\in \mathcal P\Longrightarrow \psi^+\in \mathcal P\Longrightarrow \psi\in \mathcal P 
\label{i-3}
\end{equation}
\begin{eqnarray}
\psi\in \mathcal P \ \text{ is invertible in the sense of Ehrenpreis}  \Longrightarrow \\
\Longrightarrow \psi^+\in \mathcal P\ \text{ is invertible in the sense of Ehrenpreis }
 \Longrightarrow \\ 
\Longrightarrow\psi_1\in
\mathcal P
\ \text{ is invertible in the sense of Ehrenpreis } .
\label{i-4}
\end{eqnarray}
To obtain
 (\ref{i-4}) we have to use the version (\ref{treb-2}) of the definition of slowly decreasing function 
because of the estimates of the functions $\psi^+$ and $\psi_1$ hold for non-real $z=x+2\mathrm{i} M_0\mathrm{ln}\,  |x|,$
and $z=x-2\mathrm{i} M_0\mathrm{ln}\,  |x|,$ correspondingly.

The implications (\ref{i-1}), (\ref{i-2}), (\ref{i-3}), (\ref{i-4}) lead to the required assertion. 

Q.E.D.

\smallskip

{\bf Remark 1.}
{\it Further, we will formulate and prove all assertions for entire functions with only real zeros.
By Lemma \ref{lem-1}, it is clear that
each of them has natural generalization for the case of entire functions which zeros
$\{\lambda_j\}$ are not necessarily real, but
satisfy the relation $\mathrm{Im}\, \lambda_j =O(\mathrm{ln}\, |\lambda_j|),$ $j\to\infty .$ }

\smallskip

 Let $\psi\in \mathcal P$ and $\{(a_j;m_j)\}$ denote its zero set, where
 $m_j$ denotes the multiplicity of the zero  $a_j\in\mathbb C .$ 
 In the paper \cite[Proposition 6.1]{Ehren} L. Ehrenpreis shows that  
if this function is
is  slowly decreasing function then
the inequality
\begin{equation}
\varliminf_{j\to\infty}\frac{m_j}{|\mathrm{Im}\, a_j|+\mathrm{\ln}\, |\mathrm{Re}\, a_j|} <\infty
\label{eh-nec}
\end{equation}
holds.

Given sequence  $\mathcal M=\{ \mu_k\}\subset\mathbb C$,  $|\mu_1|\le |\mu_2|\le\dots$, we denote by
 $m(z,t)$ the number of its points $\mu_k$ in the closed disc of   radius $t$ centered at 
 $z$.

Assuming that the zero set $\mathcal M\subset\mathbb R$, we improve the cited result due to L. Ehrenpreis.

\begin{lemma}
Let
$\psi\in\mathcal P$ 
be invertible in the sense of Ehrenpreis (equivalently, $\psi$ is slowly decreasing)
with zeros  $\mathcal M=\{ \mu_k\}\subset\mathbb R$.

Then, 
\begin{equation}
\varlimsup_{|x|\to\infty} \frac{m(x,1)}{\mathrm{ln}\, |x|} <\infty .
\label{star-00}
\end{equation}
\label{lem-2}
\end{lemma}

{\bf Proof.}

{\sl Wlog}, we  assume that 
 $\psi$ 
is bounded on the real axis and its exponential type equals 1.

Suppose that (\ref{star-00}) fails, that is 
 \begin{equation}
\lim_{j\to\infty} \frac{m(x_j,1)}{\mathrm{ln}\, |x_j|} =\infty 
\label{opp-1}
\end{equation}
is true for some $x_j$,  $|x_j| \to \infty .$
For clarity, assume that $x_j>0$;
and set 
\begin{equation*}
\psi_j (z)=\psi(z)(z-x_j)^{m_j}\cdot\prod\limits_{k: |\mu_k-x_j|\le 1}(z-\mu_k)^{-1}, 
\end{equation*}
where $m_j =m(x_j ,1),$ $j=1,2,\dots $
It is easy to check that  $\psi_j$ are the entire functions of the exponential type 1,
and the estimates   
$$
\sup\limits_{x\in\mathbb R} |\psi_j(x)|
\le C_0 2^{m_j},  \quad j=1,2,\dots,
$$
hold
with  $C_0=\sup\limits_{t\in\mathbb R} |\psi (t)|.$
By well-known theorem due to S. Bernstein  \cite[Chapter 11]{Boas}, the following estimates
\begin{equation}
\sup\limits_{x\in\mathbb R} |\psi^{(n)}_j (x)|\le C_0 2^{m_j} 
\label{est-deriv}
\end{equation}
are also true  for all $n, \, j\in\mathbb N.$

Further, we modify L.Ehrenpreis' argument which he used in \cite[Proposition 6.1]{Ehren}.

By the Taylor expansion of the function $\psi_j$ at
   $x_j$ and the estimates (\ref{est-deriv}), we derive that
\begin{equation*}
|\psi_j(z)|\le C_02^{m_j} (m_j!)^{-1} |z-x_j|^{m_j}e^{|z-x_j|}, \quad z\in\mathbb C.
\end{equation*}
Hence, 
for all $x\in\mathbb R$ satisfying the condition
\begin{equation}
\mathrm{ln}\, C_0 +m_j +|x-x_j|+m_j\mathrm{ln}\, |x-x_j|-\mathrm{ln}\, (m_j!)\le -l\mathrm{ln}\, x_j -m_j\mathrm{ln}\, 2 
\label{star}
\end{equation}
we have the inequality 
\begin{equation}
|\psi_j(x)|\le x_j^{-l}\cdot 2^{-m_j} ;
\label{ca-0}
\end{equation}
 here, $ l\in\mathbb N  .$

By Stirling's formula,  the relation (\ref{star}) will follow from the inequality
\begin{equation*}
|x-x_j|+m_j\mathrm{ln}\, |x-x_j|-m_j\mathrm{ln}\, m_j \le -l\mathrm{ln}\, | x_j| -C_1 m_j,
\end{equation*}
where $C_1$ is an absolute constant.

Because of  (\ref{opp-1}), for each $l\in\mathbb N$ we can find   $j_l$ such that
\begin{equation}
-l\, \mathrm{ln}\, x_j\ge -m_j , \quad j=j_l, j_l+1,\dots 
\label{cross}
\end{equation}
Let $b\in (0;1) $ and
$b<e^{-C_1-2}.$
The estimates 
(\ref{ca-0}) hold for $j\ge j_l$ and  all $x\in\mathbb R$  such that
$|x-x_j|\le b m_j .$
From this fact, inequalities  (\ref{cross}) and the relations
$$
|\psi (z)|\le 2^{m_j}|\psi_j(z) |,\quad z\in\mathbb C,\ \ j=1,2, \dots ,
$$
  we derive that
$$
|\psi (x)|\le |x_j|^{-l},\quad \text{for} \  x\in\mathbb R \quad |x-x_j|\le b l\,\mathrm{ln} x_j,\ \ j\ge j_l,\ \ l\in\mathbb N.
$$
It means that  $\psi$ is not slowly decreasing and leads to the contradiction.

Q.E.D.

\medskip

Let $\varphi\in\mathcal P$ and   
$\Lambda =\{\lambda_j\}\subset\mathbb R\setminus \{0\}$ be its zero set. 
We introduce the following notations. 
$n(z,t)$ denotes the number of points $\lambda_j$ in the closed disc of  radius $t$
centered at $z;$
$\nu (t)$ is the number of points $\lambda_j$  in  $(0;t]$ as $t> 0$,
and $(-\nu (t))$ is the number of points $\lambda_j$  in  $[t;0)$ as $t< 0;$
$2\Delta=\lim\limits_{j\to\infty}\frac{j}{|\lambda_j|}.$

\begin{theorem}
If $\varphi\in\mathcal P$ is invertible in the sense of 
Ehrenpreis and has only real zeros then
\begin{equation}
\nu (t)-\Delta t= O(\mathrm{ln}^2\, |t| ), \quad \text{as}\  |t|\to\infty .
\label{m-tm-1}
\end{equation}
\label{tm-1}
\end{theorem}

{\bf Proof.}

Below, we list some auxiliary facts.
 
F1)  Invertibility in the sense of Ehrenpreis for the function $\varphi\in\mathcal P$ with real zeros
is equivalent to the existence of constants $M_0>0$ and $r_0>1$ such that
\begin{equation}
\mathrm{ln}\, |\varphi (z)|\ge-M_0 \mathrm{\ln}\, |x|,\quad z=x+\mathrm{i}y, \quad |x|\ge r_0, \quad |y|\ge M_0 \mathrm{ln}\, |x| .
\label{F1}
\end{equation}

F2) \cite[Section 3]{Ber-Tayl}.  If  $\varphi\in\mathcal P$ is  invertible 
in the sense  of Ehrenpreis and all its zeros are real then the set
\begin{equation}
\{z: \ \mathrm{ln}\, |\varphi (z)|<- M_0 \mathrm{ln}\, |z|, \
|x|\ge r_0, \quad |y|\le M_0 \mathrm{ln}\, |x| \}
\label{ex-set}
\end{equation}
consists of relatively compact components $G$ of the diameter
 \begin{equation}
d_G\le M\mathrm{ln}\, |z|,\ \forall z\in G.
\label{dia}
\end{equation}

F3)  By Lemma \ref{lem-2} and (\ref{dia}),   
the number of points  $\lambda_j \in G$ is 
 $O(\mathrm{ln}^2\, |z |)$
$\forall z\in G ,$
for each relatively compact component $G$ of the set (\ref{ex-set}).

F4) Let  $\varphi\in\mathcal P$ be invertible in the sense of Ehrenpreis and have only real zeros.
By standard technique of the theory of entire functions,
 from F3), we derive the estimate
\begin{equation}
\mathrm{ln}\, |\varphi (z)|\ge -C\mathrm{ln}^2\, |z|\mathrm{ln}\,\mathrm{ln}\, |z| 
\label{F4}
\end{equation}
if $z=x+\mathrm{i} y,$ $ |x|\ge r_0,$ $ |y|\ge r_0$.
Here,  $C>0$ depends only on  $M_0,$ $r_0 $ and the density  of $\Lambda$.

F5) Let  $\varphi\in\mathcal P $ be such as above.
By Theorem III.G.1  \cite{Koosis-1}, 
\begin{equation}
\mathrm{ln}\, |\varphi (z)| =\pi\Delta\mathrm{Im} z+\frac{1}{\pi}\int_{-\infty}^{\infty} \frac{\mathrm{Im} z \mathrm{ln}\, |\varphi (t)|}{|z-t|^2}
\mathrm{d} t, \quad \mathrm{Im} z >0,
\label{r-1}
\end{equation}
\begin{equation}
\mathrm{ln}\, |\varphi (z)| =-\pi\Delta\mathrm{Im} z-\frac{1}{\pi}\int_{-\infty}^{\infty} \frac{\mathrm{Im} z \mathrm{ln}\, |\varphi (t)|}{|z-t|^2}
\mathrm{d} t, \quad \mathrm{Im} z <0.
\label{r-2}
\end{equation}

\medskip

Now, we estimate the difference $(\nu (x_0)-\Delta x_0)$ 
for  an arbitrary  $x_0>r_0,$ $x_0\not\in\Lambda .$
According to the well-known formula, we have
\begin{equation}
\nu (x_0) =\frac{1}{2\pi\mathrm{i}} \int_{\Gamma^+} \frac{\varphi' (z)}{\varphi (z)}\mathrm{d} z,
\label{p-1}
\end{equation}
where   $\Gamma^+$ is the boundary of the rectangle
$$
\{z=x+\mathrm{i}y :\ 0\le x\le x_0,\ |y|\le 3M_0\mathrm{ln}\, x_0 \}.
$$
 
 Taking into account the properties of $\varphi,$
for the real parts of both sides of (\ref{p-1})
we obtain
\begin{multline*}
\nu (x_0)=\frac{1}{2\pi}\bigl(
\mathrm{arg}\, \varphi (x_0-3\mathrm{i}M_0\mathrm{ln}\, x_0)
-\mathrm{arg}\, \varphi (-3\mathrm{i}M_0\mathrm{ln}\, x_0)+\\
+\mathrm{arg}\, \varphi (3\mathrm{i}M_0\mathrm{ln}\, x_0)
-\mathrm{arg}\, \varphi (x_0+3\mathrm{i}M_0\mathrm{ln}\, x_0)
\bigr) .
\end{multline*}
Here, $\mathrm{arg}\, \varphi$ is the imaginary part of some branch
(not necessarily, the main one)
of the function
$\left(\mathrm{ln} \, \varphi\right)$ which is analytic in   each of the 
 half-planes  $\mathrm{Im}\, z<0$ and $\mathrm{Im}\, z>0$.
Generally speaking, we choose different branches
for each half-plane.

By  (\ref{r-1}) and (\ref{r-2}), we see that
$$
\mathrm{ln}\, |\varphi (z)| =\mathrm{Re}\, \left(-\mathrm{i}\pi\Delta z
+\frac{\mathrm{i}}{\pi}\int_{-\infty}^{\infty}\left( \frac{1}{z-t}+\frac{t}{t^2+1}\right)\mathrm{ln}\, |\varphi (t)|\mathrm{d}t\right),
\quad \mathrm{Im}\, z>0,
$$
$$
\mathrm{ln}\, |\varphi (z)| =\mathrm{Re}\, \left(\mathrm{i}\pi\Delta z
-\frac{\mathrm{i}}{\pi}\int_{-\infty}^{\infty}\left( \frac{1}{z-t}+\frac{t}{t^2+1}\right)\mathrm{ln}\, |\varphi (t)|\mathrm{d}t\right),
\quad \mathrm{Im}\, z<0 .
$$
Hence,
$$
\mathrm{arg}\, |\varphi (z)| =\mathrm{Im}\, \left(-\mathrm{i}\pi\Delta z
+\frac{\mathrm{i}}{\pi}\int_{-\infty}^{\infty}\left( \frac{1}{z-t}+\frac{t}{t^2+1}\right)\mathrm{ln}\, |\varphi (t)|\mathrm{d}t\right) +\mathrm{const},
$$
  as $\mathrm{Im}\, z>0;$
$$
\mathrm{arg}\, |\varphi (z)| =\mathrm{Im}\, \left(\mathrm{i}\pi\Delta z
-\frac{\mathrm{i}}{\pi}\int_{-\infty}^{\infty}\left( \frac{1}{z-t}+\frac{t}{t^2+1}\right)\mathrm{ln}\, |\varphi (t)|\mathrm{d}t\right) +\mathrm{const},
 $$
as $\mathrm{Im}\, z<0$
(\cite[III.H.2]{Koosis-1}).

Generally speaking, there are different constants
in the right-hand sides of two last formulas.

From the above, taking into account the relation
 $\overline{\varphi (\bar z)}=\varphi (z),$
we conclude that
\begin{multline}
\nu (x_0) =-\frac{1}{\pi} \left(\mathrm{arg}\, \varphi (x_0+3\mathrm{i}M_0\mathrm{ln}\, x_0)
-\mathrm{arg}\, \varphi (3\mathrm{i}M_0\mathrm{ln}\, x_0)\right) +\mathrm{const}=\\
=\Delta x_0-\int_{-\infty}^{\infty} \left( \frac{x_0-t}{(x_0-t)^2+9M_0^2\mathrm{ln}^2\, x_0}+\frac{t}{t^2+1}\right)
\mathrm{ln}\, |\varphi (t)|\mathrm{d}t+\\
\int_{-\infty}^{\infty} \left( \frac{-t}{t^2+9M_0^2\mathrm{ln}^2\, x_0}+\frac{t}{t^2+1}\right)
\mathrm{ln}\, |\varphi (t)|\mathrm{d}t .
\label{r-3}
\end{multline}
Last formula does not completely suit for estimating  the expression $(\nu (x_0)-\Delta x_0)$,
because there is no finite estimates from below for the function $\mathrm{ln}\, |\varphi (t)| $ 
on the whole real line.
To workaround this issue, we consider the function
$\psi (\tilde{z})=\varphi (z),$
where
$\tilde{z}=z-2\mathrm{i}M_0\mathrm{ln}\,x_0.$
 The function $\psi$ is analytic and non-vanishing 
in the {\sl closed} upper half-plane
$\mathrm{Im}\, \tilde{z}\ge 0.$
And we have the representation similar to  (\ref{r-1})
for $\mathrm{ln}\, |\psi (\tilde{z})|:$ 

\begin{multline*} 
\mathrm{ln}\, |\psi (\tilde{z})| =\pi\Delta \mathrm{Im}\, \tilde{z} +\frac{1}{\pi}\int_{-\infty}^{\infty}
\frac{\mathrm{Im}\,\tilde{z}\, \mathrm{ln}\, |\psi (t)|}{|\tilde{z}-t|^2}\mathrm{d}t =\\
= \mathrm{Re}\,\left( -\mathrm{i}\pi\Delta \tilde{z} +\frac{\mathrm{i}}{\pi}
\int_{-\infty}^{\infty} \left( \frac{1}{\tilde{z}-t}+\frac{t}{t^2+1}\right)\mathrm{ln}\, |\psi (t)|\mathrm{d} t
\right), \qquad \mathrm{Im}\, \tilde{z} > 0.
\end{multline*}
Now, we rewrite this representation in terms of  $\varphi$ and $z$:
\begin{multline*}
\mathrm{ln}\, |\varphi (z)| =  \mathrm{Re}\,\Biggl( -\mathrm{i}\pi\Delta (z -2\mathrm{i} M_0\mathrm{ln}\, x_0)+\\
+\frac{\mathrm{i}}{\pi}
\int_{-\infty}^{\infty} \Bigl( \frac{1}{z-2\mathrm{i}M_0\mathrm{ln}\, x_0-t}+\frac{t}{t^2+1}\Bigr)
\mathrm{ln}\, |\varphi (t+2\mathrm{i}M_0\mathrm{ln}\, x_0)|\mathrm{d} t
\Biggr) 
\end{multline*}
if $ \mathrm{Im}\, z > 2M_0\mathrm{ln}\, x_0 .$
By the same way as it has been done for (\ref{r-3}),  we  obtain 
\begin{multline}
\nu (x_0) =\Delta x_0 -\int_{-\infty}^{\infty}
\left( \frac{x_0-t}{(x_0-t)^2+M_0^2\mathrm{ln}^2\, x_0}+\frac{t}{t^2+1}\right)
\mathrm{ln}\, |\varphi (t+2\mathrm{i} M_0\mathrm{ln}\, x_0)|\mathrm{d}t +\\
+\int_{-\infty}^{\infty}
\left( \frac{-t}{t^2+M_0^2\mathrm{ln}^2\, x_0}+\frac{t}{t^2+1}\right)
\mathrm{ln}\, |\varphi (2\mathrm{i} M_0\mathrm{ln}\, x_0)|\mathrm{d}t .
\label{r-4}
\end{multline} 
It follows that
\begin{equation}
|\nu_0(x_0)-\Delta x_0| \le |I_1|+|I_2|,
\label{e-e}
\end{equation}
where
$$
I_1=\int_{-\infty}^{\infty}
\left( \frac{x_0-t}{(x_0-t)^2+M_0^2\mathrm{ln}^2\, x_0}+\frac{t}{t^2+1}\right)
\mathrm{ln}\, |\varphi (t+2\mathrm{i} M_0\mathrm{ln}\, x_0)|\mathrm{d}t ,
$$
$$
I_2=\int_{-\infty}^{\infty}
\left( \frac{-t}{t^2+M_0^2\mathrm{ln}^2\, x_0}+\frac{t}{t^2+1}\right)
\mathrm{ln}\, |\varphi (2\mathrm{i} M_0\mathrm{ln}\, x_0)|\mathrm{d}t .
$$
{\sl Wlog}, we assume that
$$
\mathrm{ln}\, |\varphi (z)|\le \pi\Delta |\mathrm{Im}\, z|, \quad z\in\mathbb C.
$$
Taking into account (\ref{F1}), we have
\begin{equation}
\left|\mathrm{ln} \, |\varphi (z)| \right|\le C_1\mathrm{ln}\, x_0,
\label{es-1}
\end{equation}
for all $z=t+2\mathrm{i}M_0\mathrm{ln}\, x_0$ with $|t|\le x_0^2$;
here, the positive constant $C_1$ depends only on $r_0,$ $M_0,$ $\Delta.$

Further, by (\ref{F4}), we derive that 
\begin{equation}
\left|\mathrm{ln} \, |\varphi (z)| \right|\le C_2\mathrm{ln}^2\, |t| \mathrm{ln}\,\mathrm{ln}\, |t|
\label{es-2}
\end{equation}
if $z=t+2\mathrm{i}M_0\mathrm{ln}\, x_0$
and $ |t|\ge x_0^2$,
where the constant $C_2>0$
depends only on $r_0,$ $M_0,$ $\Delta.$

Now,  we estimate $|I_1|$
by help of (\ref{es-1})
and (\ref{es-2}):
\begin{multline}
|I_1|\le \left|  
\int_{|t|\le x_0^2}
\left( \frac{x_0-t}{(x_0-t)^2+M_0^2\mathrm{ln}^2\, x_0}+\frac{t}{t^2+1}\right)
\mathrm{ln}\, |\varphi (t+2\mathrm{i} M_0\mathrm{ln}\, x_0)|\mathrm{d}t 
\right|+\\
+\left|\int_{|t|\ge x_0^2}
\left( \frac{x_0-t}{(x_0-t)^2+M_0^2\mathrm{ln}^2\, x_0}+\frac{t}{t^2+1}\right)
\mathrm{ln}\, |\varphi (t+2\mathrm{i} M_0\mathrm{ln}\, x_0)|\mathrm{d}t \right|=\\
=|I_{11}|+|I_{12}|.
\label{es-3}
\end{multline}
Taking into account (\ref{es-1}), we get
\begin{multline*}
|I_{11}|\le C_1\mathrm{ln}\, x_0\left( \int_{-x_0^2}^{x_0^2}
\frac{|x_0-t|}{(x_0-t)^2+M_0^2\mathrm{ln}^2\, x_0 }
\mathrm{d}t
 +\int_{-x_0^2}^{x_0^2} \frac{|t|}{t^2+1}
\mathrm{d}t
\right)\le \\ \le \widetilde{C}_1 \mathrm{ln}^2\, x_0,
\end{multline*}
where  $\widetilde{C}_1>0$ depends only on $r_0,$ $M_0,$ $\Delta. $

By simple transformations, we see that
$$\left| \frac{x_0-t}{(x_0-t)^2+M_0^2\mathrm{ln}^2\, x_0}+\frac{t}{t^2+1}\right|
\le \frac{C_3}{|t|^{3/2}}, \quad |t|\ge x_0^2,
$$
where $C_3$ depends only on $r_0,$ $M_0,$ $\Delta .$

 (\ref{es-2}) means that
$$
\left|\mathrm{ln}\, |\varphi (t+2\mathrm{i} M_0\mathrm{ln}\, x_0)|\right|
\le C_2\mathrm{ln}^2\, |t| \mathrm{ln}\,\mathrm{ln}\, |t|.
$$
Summarizing the above, we obtain that
$$
|I_{12}|\le C_4,
$$
where  $C_4>0$ 
depends only on $r_0,$ $M_0,$ $\Delta. $

From the representation (\ref{es-3}) and the estimates for  $I_{11}$ and $I_{12}$ 
it follows that
$$
I_1 =O(\mathrm{ln}^2\, x_0^2),\quad x_0\to+\infty .
$$
 Similarly, we get
$$
I_2 =O(\mathrm{ln}^2\, x_0^2),\quad x_0\to+\infty .
$$
Now, the inequality (\ref{e-e}) and  two last estimates lead to the required relation
\begin{equation}
\nu(x)-\Delta x =O(\mathrm{ln}^2\, x), \quad x\to+\infty , \quad x\not\in \Lambda .
\label{star-1}
\end{equation}

By the same way, we manage with the case $x<0.$

At last, for $x\in\Lambda$, the asymptotic relation (\ref{m-tm-1}) follows from (\ref{star-1}) and Lemma \ref{lem-2}.
Q.E.D.
\smallskip

{\bf Remark 2. } 
{\it Set 
$l(t)=\mathrm{ln}\, (1+t^2),$ $t\in\mathbb R,$
$\lambda_j=j+l(|j|),$ $j=\pm 1,$ $\pm 2,$ $\dots $
It is known that the function
$$
\varphi (z)=\lim_{R\to+\infty}\prod\limits_{|\lambda_j|< R}\left( 1-\frac{z}{\lambda_j}\right)
$$ 
is invertible in the sense of Ehrenpreis in the algebra $\mathcal P$ \cite[Theorem 1]{NF-POMI};
and, at the same time, 
$$
\nu (x) =[x- \mathrm{ln}\, (1+x^2)+o(1) ],\quad x\to+\infty 
 $$
(see \cite[Lemma 1]{NF-POMI}). 
Hence, the necessary condition  (\ref{m-tm-1}) cannot be improved.
}

\smallskip

{\bf Remark 3.}
{\it Generally speaking, the necessary condition of invertibility in the sense of Ehrenpreis (\ref{m-tm-1}) 
is not the sufficient one.
Let us demonstrate it.
Set
$$
\varphi_0(z)=\frac{\sin \pi z}{\pi z s_0(z)},
$$
where $s_0(z)=\prod\limits_{k=1}^{\infty}\left(1-\frac{z^2}{4^k}\right)$.
By well-known estimates for the functions $\sin\pi z$
and $s_0$  \cite[Chapter 3]{Boas},
we see that 
$$
\mathrm{ln}\, |\varphi_0 (x)|\le -C\mathrm{ln}^2\, (2+|x|),\quad x\in\mathbb R,
$$
where $C>0$.
Hence, $\varphi_0$ is not slowly decreasing function; equivalently, it is not invertible in the sense of Ehrenpreis.
Nevertheless, the corresponding function $\nu$ 
satisfies the condition
$$
\nu (x)-x= O(\mathrm{ln}\, |x|) \quad  \text{as}\ |x|\to\infty.
 $$ 
}

\section{Criteria of  invertibility in the sense of Ehrenpreis}

\subsection{The zero set is an even sequence.}

We will use the following lemma due to S. Favorov
 \cite[Lemma 1]{Fav}.

\smallskip
{\bf Lemma B.}
{\it Let  $A=\{ a_j\}\subset \mathbb C\setminus\{ 0\}$ be the sequence satisfying the conditions
\begin{eqnarray}
\exists \lim_{R\to\infty }\sum_{|a_j|<R} a_j^{-1},
\label{f1}\\
n_A(0,t) =O(t),\quad t\to\infty, 
\label{f2}\\
n_A(0,t+1)-n_A(0,t) =o(t), \quad t\to\infty,
\label{f3}
\end{eqnarray}
where $n_A(z,t)$ denotes the number of points $a_j$ in the closed disc of  radius $t$
centered at $z.$ 
Then,  
\begin{equation}
g (z)=\lim_{R\to\infty }\prod\limits_{|a_j|\le R}\left( 1-\frac{z}{a_j}\right)
\label{def-g}
\end{equation}
is the entire function of exponential type,
and 
\begin{equation}
\mathrm{ln}\, |g (z)| =\int_0^{\infty} \frac{n_A(0,t)-n_A(z,t)}{t}\mathrm{d} t
\label{log-g}
\end{equation} 
holds for all $z\in\mathbb C$. }
\smallskip

For $\Lambda^+ =\{\lambda_j\},$ $0<\lambda_1\le\lambda_2\le\dots ,$
such that
$$
\exists \lim_{j\to\infty} \frac{j}{\lambda_j}=\Delta,
$$
the corresponding even sequence $\Lambda =\Lambda^+\bigcup (-\Lambda^+)$
satisfies the assumptions of Lemma B. 
Hence, the formula
\begin{equation}
\varphi (z)=\prod\limits_{j=1}^{\infty}\left( 1-\frac{z^2}{\lambda^2_j}\right)
\label{def-fi}
\end{equation}
defines  an entire function of exponential type  $\pi \Delta .$ 

By help of Theorem \ref{tm-1},  we obtain the necessary and sufficient conditions
under which $\varphi $ is invertible in the sense of Ehrenpreis element of the algebra $\mathcal P.$
To arrive to this result, first, we prove some auxiliary assertions.

For $x\in\mathbb R$ and $t>0$,  we denote by $n^+(x;t)$ and $n^-(x;t)$ 
the numbers of points of $\Lambda$ laying in the intervals 
 $(x;x+t]$ and $(x-t;x],$ correspondingly.

\begin{lemma}
If
$$
n^+(0;x)-\Delta x =O(\mathrm{ln}^2\, x) \ \text{as}\ x\to +\infty
$$
then 
\begin{equation}
\int_{x\mathrm{ln}\, x}^{+\infty} \frac{n^+(t;x) -n^-(t;x)}{t}\mathrm{d} t =O(\mathrm{ln}\, x) \ \text{as}\  x\to+\infty .
\label{int-1}
\end{equation}
\label{lem-3}
\end{lemma} 

{\bf  Proof.}
We fix an arbitrary  $\sigma >2$ and write the representation 
\begin{multline}
\int_{x\mathrm{ln}\, x}^{+\infty} \frac{n^+(t;x) -n^-(t;x)}{t}\mathrm{d}\, t =
\int_{x\mathrm{ln}\, x}^{x^\sigma} \frac{n^+(t;x) -n^-(t;x)}{t}\mathrm{d}\, t +\\
+\int_{x^{\sigma}}^{+\infty} \frac{n^+(t;x) -n^-(t;x)}{t}\mathrm{d}\, t =J_1+J_2.
\label{rep}
\end{multline}
 By  the evenness of $\Lambda $ and Theorem \ref{tm-1}, we get
\begin{multline*}
J_2=\int_{x^{\sigma}}^{+\infty} \frac{xn^+(t;x)}{t(t+x)}\mathrm{d}\, t-\int_{x^{\sigma}-x}^{x^{\sigma}}
\frac{n^{+}(t;x)}{t+x}\mathrm{d} t
=\\
 =\int_{x^{\sigma}}^{+\infty}
  \frac{O(x^2)+xO(\mathrm{ln}^2\, t)}{t(t+x)}\mathrm{d}\, t 
+O\left(x\mathrm{ln}\, \left( 1+\frac{x}{x^{\sigma}-x}\right)\right)=
\\
=O(1), \quad x\to+\infty.
\end{multline*}

Taking into account   the evenness of $\Lambda $ and Theorem \ref{tm-1} one more time, we estimate  $J_1:$
\begin{multline*}
\int_{x\mathrm{ln}\, x}^{x^\sigma} \frac{n^+(t;x) -n^-(t;x)}{t}\mathrm{d} t =
-\int_{x\mathrm{ln}\, x-x}^{x\mathrm{ln}\, x} \frac{n^+(t;x)\pm\Delta x}{t+x}\mathrm{d} t +\\
\int_{x\mathrm{ln}\, x}^{x^{\sigma}- x} n^+(t;x)\left(\frac{1}{t} -\frac{1}{t+x}\right)\mathrm{d} t
+\int_{x^{\sigma}-x}^{x^{\sigma}} \frac{n^+(t;x)}{t}\mathrm{d}\, t=\\
= O(\mathrm{ln}\, x)-\Delta x\mathrm{ln}\, \left( 1+\frac{1}{\mathrm{ln}\, x}\right)+
\int_{x\mathrm{ln}\, x}^{x^{\sigma}- x} (O(\mathrm{ln}^2\, t)+\Delta x)\left(\frac{1}{t} -\frac{1}{t+x}\right)\mathrm{d} t+\\
+\int_{x^{\sigma}-x}^{x^{\sigma}}(O(\mathrm{ln}^2\, t)+\Delta x)\left(\frac{1}{t} -\frac{1}{t+x}\right)\mathrm{d} t=\\
=O(\mathrm{ln}\, x)-\Delta x\mathrm{ln}\, \left( 1+\frac{1}{\mathrm{ln}\, x}\right)+O(\mathrm{ln}\, x)+\Delta x
\mathrm{ln}\, \frac{x^{\sigma}-x}{x}-\\
-\Delta x \mathrm{ln}\, \left( 1-\frac{1}{\mathrm{ln}\, x+1}\right) +O(1)+\Delta x \mathrm{ln}\, \frac{x^{\sigma}}{x^{\sigma}-x}
=O (\mathrm{ln}\, x), \quad x\to+\infty .
\end{multline*}

From (\ref{rep}) and the estimates for $J_1$ and $J_2$, it follows that  (\ref{int-1}) holds.

Q.E.D.

\begin{lemma}
Let $\Lambda =\Lambda^+\bigcup\Lambda^-,$ where $\Lambda^+=\{\lambda_j\}$ is positive sequence 
of density $\Delta .$

Then, for a fixed $A>0$, the relation 
\begin{equation}
I:=\int_{|x|\mathrm{ln}\, |x|}^{\infty} \frac{n(x,t)-n(x+\mathrm{i}A\mathrm{ln}\, |x|;t)}{t}\mathrm{d}t =O(A^2), \quad |x|\to\infty ,
\label{l-4}
\end{equation}
holds;
here  $n(z,t)$ denotes the number of points $\pm\lambda_j$ in the closed disc of 
radius $ t$ centered at $z.$  
\label{lem-4}
\end{lemma}

{\bf Proof.}

Assume that $x>0$ (for  $x<0$ the same argument works).

We have
$$
0\le I=\int_{|x|\mathrm{ln}\, |x|}^{\infty} \frac{n^+(t+x-r_t;r_t)}{t}\mathrm{d} t+
\int_{|x|\mathrm{ln}\, |x|}^{\infty} \frac{n^+(-t+x;r_t)}{t}\mathrm{d} t =I_1+I_2,
$$ 
where $r_t=\sqrt{t^2 -A^2\mathrm{ln}^2\, x}$.

First, we estimate  $I_1:$ 
\begin{multline*}
0\le I_1=\sum_{\lambda_j\ge x+x\mathrm{ln}\, x} \int_{\lambda_j-x}^{\sqrt{\lambda_j^2+A^2\mathrm{ln}^2\, x }-x}
\frac{\mathrm{d}t}{t}\le
\\ 
\le
\mathrm{const}\, 
\sum_{ \lambda_j\ge x+x\mathrm{ln}\, x} \frac{A^2\mathrm{ln}^2\, x}{(\lambda_j-x)^2}\le C_0 A^2, \quad x\ge x_0>0;
\end{multline*}
here, the positive constant $C_0$ depends only on  $\{\lambda_j\}$.

Second integral $I_2$ is estimated by the similar way.
 
From the estimates for  $I_1,$ $I_2$, we obtain the required one for $I=I_1+I_2.$

Q.E.D.
\medskip

Now, we are ready to formulate and prove the main result.

\begin{theorem}
Let $\Lambda =\{ \pm\lambda_j\},$ $0<\lambda_1\le\lambda_2\le\dots,$ 
and the finite limit exists $\lim\limits_{j\to\infty} \frac{j}{\lambda_j}=:\Delta .$

The entire function
\begin{equation}
\varphi (z)=\prod_{j=1}^{\infty}\left( 1-\frac{z^2}{\lambda_j^2}\right) 
\label{d-fi}
\end{equation}
belongs to the algebra $\mathcal P$ and is invertible in the sense of 
Ehrenpreis if and only if the following two conditions hold:
\begin{equation}
n^{+} (0;x)-\Delta x =O(\mathrm{ln}^2\, x),\quad x\to\infty,
\label{cond-1}
\end{equation}
\begin{equation}
\varlimsup_{A\to\infty}\varlimsup_{x\to\infty} \frac{1}{A\mathrm{ln}\, x}\left|
\int_{A\mathrm{ln}\, x}^{x\mathrm{ln}\, x} \frac{n(0,t)-n(x+\mathrm{i}A\mathrm{ln}\, x ,t)}{t}\mathrm{d} t
\right| <+\infty .
\label{cond-2}
\end{equation}
\label{tm-2}
\end{theorem}

{\bf Proof.}

{\sl Necessity.} 

Assume that the  function
 $\varphi$ defined by (\ref{d-fi}) is  invertible in the sense of Ehrenpreis element
of the algebra $\mathcal P$. 
From Theorem \ref{tm-1}, it follows that the relation (\ref{cond-1}) is true.

{\sl Wlog}, we may assume that
$|\varphi (x)|\le 1$ and 
$$
\mathrm{ln}\, |\varphi (z)|\le \pi\Delta |\mathrm{Im}\, z|, \quad z\in\mathbb C.
$$
With Lemma B, it leads to the inequality 
 $$
\int_0^{\infty} \frac{n(0,t)-n(x +\mathrm{i} A\mathrm{ln}\, x,t)}{t}\mathrm{d} t\le\pi\Delta A\mathrm{ln}\, x 
$$
for all $A>0$ and $x\ge x_0>0.$ 

Hence,  we have
\begin{multline*}
\int_{A\mathrm{ln}\, x}^{x\mathrm{ln}\, x}\frac{n(0,t)-n(x +\mathrm{i} A\mathrm{ln}\, x,t)}{t}\mathrm{d} t\le \\
\le \pi\Delta A\mathrm{ln}\, x+\int_{x\mathrm{ln}\, x}^{\infty} \frac{n(x +\mathrm{i} A\mathrm{ln}\, x,t)-n(0,t)}{t}\mathrm{d} t\le \\
\le \pi\Delta A\mathrm{ln}\, x+\int_{x\mathrm{ln}\, x}^{\infty} \frac{n(x,t)-n(0,t)}{t}\mathrm{d} t\le
\pi\Delta A\mathrm{ln}\, x+\mathrm{const }\, \mathrm{ln}\, x
\end{multline*}
for all $x\ge x_0;$
and in this chain of inequalities, the last one is true because of Lemma \ref{lem-3}.
The above estimates give us the following inequality
\begin{equation}
\varlimsup_{A\to\infty}\varlimsup_{x\to\infty} \frac{1}{A\mathrm{ln}\, x}
\int_{A\mathrm{ln}\, x}^{x\mathrm{ln}\, x} \frac{n(0,t)-n(x+\mathrm{i}A\mathrm{ln}\, x ,t)}{t}\mathrm{d} t
\le \pi\Delta .
\label{cond-2-1}
\end{equation}

Further, there exists  $A_0>0$ such that
$$
\mathrm{ln}\, |\varphi (x+\mathrm{i} A_0\mathrm{ln}\, x)|\ge -A_0\mathrm{ln}\, x, \quad x\ge x_0.
$$
All zeros of  $\varphi$ are real. It implies  
the estimate
$$
\mathrm{ln}\, |\varphi (x+\mathrm{i} A\mathrm{ln}\, x)|\ge -A_0\mathrm{ln}\, x, \quad A\ge A_0,\  x\ge x_0.
$$
This estimate, together with  Lemma B, Lemma \ref{lem-3}, Lemma \ref{lem-4} lead to
\begin{multline*}
\int_{A\mathrm{ln}\, x}^{x\mathrm{ln}\, x} \frac{n(x+\mathrm{i}A\mathrm{ln}\, x, t)-n(0,t)}{t}\mathrm{d} t\le
A_0\mathrm{ln}\, x+\mathrm{const}\, A\mathrm{ln}\, x+ \\
+ \int_{x\mathrm{ln}\, x}^{\infty} \frac{n(0,t)-n(x+\mathrm{i}A\mathrm{ln}\, x,t)}{t}\mathrm{d} t\le
\mathrm{const}\, \mathrm{ln}\, x+\mathrm{const}\, A\mathrm{ln}\, x +\mathrm{const}\, A^2.
\end{multline*}
Finally, 
\begin{equation}
\varlimsup_{A\to\infty}\varlimsup_{x\to\infty} \frac{1}{A\mathrm{ln}\, x}
\int_{A\mathrm{ln}\, x}^{x\mathrm{ln}\, x} \frac{n(x+\mathrm{i}A\mathrm{ln}\, x,t)-n(0 ,t)}{t}\mathrm{d} t
\le \pi\Delta .
\label{cond-2-2}
\end{equation}

The inequalities (\ref{cond-2-1}), (\ref{cond-2-2}) imply (\ref{cond-2}).

\smallskip

{\sl Sufficiency.}

From the relation  (\ref{cond-1})
it follows that $\varphi$ (defined by (\ref{d-fi})) belongs to the class C (Cartwright class of entire functions).

The condition (\ref{cond-2}), Lemma B and Lemma \ref{lem-3} give us the estimate
\begin{multline*}
\mathrm{ln}\, |\varphi (x)| = \int_{0}^{A_0\mathrm{ln}\, x} \frac{n(0,t)-n(x,t)}{t}\mathrm{d} t+\\
+
\int_{A_0\mathrm{ln}\, x}^{x\mathrm{ln}\, x} \frac{n(0,t)-n(x,t)}{t}\mathrm{d} t+\int_{x\mathrm{ln}\, x}^{\infty} 
\frac{n(0,t)-n(x,t)}{t}\mathrm{d} t \le \\
\le \mathrm{const}\, \mathrm{ln}\, x + \int_{A_0\mathrm{ln}\, x}^{x\mathrm{ln}\, x} 
\frac{n(0,t)-n(x+\mathrm{i}A_0\mathrm{ln}\, x,t)}{t}\mathrm{d} t \le 
 \mathrm{const}\, \mathrm{ln}\, x   ,\ \ x\ge x_0,
\end{multline*} 
 where $A_0>0,$ $x_0 >1.$
Hence, $\varphi \in\mathcal P.$

To prove that $\varphi$
is invertible in the sense of Ehrenpreis,
we again use Lemma B, Lemma \ref{lem-3} 
and the relation (\ref{cond-2}). 
It allows us to derive the following estimate
\begin{multline*}
\mathrm{ln}\, |\varphi (x+\mathrm{i} A_0\mathrm{ln}\, x)| \ge  
\int_{A_0\mathrm{ln}\, x}^{x\mathrm{ln}\, x} \frac{n(0,t)-n(x+\mathrm{i} A_0\mathrm{ln}\, x,t)}{t}\mathrm{d} t+\\
+\int_{x\mathrm{ln}\, x}^{\infty} \frac{n(0,t)-n(x+\mathrm{i} A_0\mathrm{ln}\, x,t)}{t}\mathrm{d} t\ge \mathrm{const}\, \mathrm{ln}\, x, \qquad x\ge x_0,
\end{multline*}
 where $A_0>0,$ $x_0 >1.$ 
This estimate means that   $\varphi$ is invertible in the sense of Ehrenpres function.

Q.E.D.

\subsection{Criterion of invertibility
in the sense of Ehrenpreis for an arbitrary function in $\mathcal P$ which zeros are real.}

Assume that $\psi\in\mathcal P$ and its zero set $\Lambda\subset \mathbb R\setminus\{ 0\}$ has density $2\Delta.$
Set
$$
L(t)=\nu (t)-\Delta t,\quad  L^*(t)=L(t)-L(-t),\quad t\in\mathbb R;
$$
here, as above, 
$\nu (t)$ denotes the number of points $\lambda_j\in\Lambda$ in the interval $(0;t]$ as $t>0,$
and $(-\nu (t))$ is  the number of points $\lambda_j\in\Lambda$ in the interval $[-t;0)$ as $t<0.$

\begin{theorem}
The function $\psi\in\mathcal P$ with zero set $\Lambda\subset \mathbb R\setminus\{ 0\}$ 
of density $2\Delta$ is invertible in the sense of Ehrenpreis
if and only if the  relations 
\begin{equation}
L(x) =O(\mathrm{ln}^2\, |x|),\quad |x|\to\infty,
\label{cond-1-3}
\end{equation}
\begin{equation}
\varlimsup_{A\to\infty}\varlimsup_{x\to\infty} \frac{1}{A\mathrm{ln}\, x}\left|
\int_{A\mathrm{ln}\, x}^{x\mathrm{ln}\, x} \frac{2L^*(t)-L^*(x+r_{t,A})+
L^*(x-r_{t,A})}{t}\mathrm{d} t
\right| <+\infty ,
\label{cond-2-3}
\end{equation}
hold.
Here, $r_{t,A}=\sqrt{t^2-A^2\mathrm{ln}^2\, x}$.
\label{tm-3}
\end{theorem}

{\bf Proof.}

First of all, we notice that because of
Theorem 2.2  \cite{Ehren}, either both functions, $\psi(z)\in\mathcal P$ 
and   $\varphi (z) : =\psi(z)\psi (-z)$  are invertible in the sense of Ehrenpreis or not.

Further, it is not difficult to see that  (\ref{cond-1-3}) is valid if and only if  (\ref{cond-1}) hold for both functions, 
$\psi$ and $\varphi .$ 
And the validity of (\ref{cond-2-3}) for $\psi$ is equivalent to the validity
of (\ref{cond-2}) for  even  function  $\varphi$.

From the above, we can conclude that  (\ref{cond-1-3}) and
(\ref{cond-2-3}) imply  invertibility in the sense
of Ehrenpreis of the function $\varphi (z) =\psi(z)\psi (-z)$ and, consequently, the same property
for $\psi .$

Conversely, if $\psi\in\mathcal P$    is invertible in the sense of Ehrenpreis
then the same is true for the function $\varphi (z) =\psi(z)\psi (-z).$
The relation  (\ref{cond-1-3}) holds because of Theorem \ref{tm-1}.
  
	Applying Theorem  \ref{tm-2} to the function  $\varphi$,
we see that
(\ref{cond-2}) is true for the counting function of its zeros.
It implies (\ref{cond-2-3}). 

Q.E.D.

\smallskip

{\bf Remark 4.}
{\it It is well-known that 
entire function
$\varphi $ belongs to the algebra $\mathcal P$ and it is slowly decreasing
if its zero set $\mathcal M=\{\mu_k\},$ $k\in\mathbb Z$, formed by bounded perturbations of the sequence of integers:
 \begin{equation}
|\mu_k-k|\le L, \quad k\in\mathbb Z,
\label{PW}
\end{equation}
for some $L>0$
(\cite[Theorem XXXIII]{PW}). 

We can derive  more general and convenient 
sufficient conditions of  invertibility in the sense of Ehrenpreis 
from Theorems \ref{tm-2} and \ref{tm-3}.

In particular, assume that for a sequence $\Lambda =\{\lambda_j\},$ $0<\lambda_1\le \lambda_2\le\dots$, 
 the following relation
\begin{equation}
\nu (t)-\Delta t=O(1)
\label{PW-1}
\end{equation}
holds as $t\to\infty .$
Evidently, it is  less  restricted  than  (\ref{PW}).
At the same time,
we can easily check that the conditions 
of the criterion in Theorem \ref{tm-2}
are satisfied for the entire function
$$
\varphi (z)=\prod_{k=1}^{\infty}\left( 1-\frac{z^2}{\lambda_j^2}\right).
$$
Hence, $\varphi \in\mathcal P $ and it  
is invertible in the sense of Ehrenpreis.

For example, both sequences

\noindent
$
1) \ \Lambda =\{\pm\lambda_j\}\bigcup\{ \pm e^{\sqrt{j}}\}_{j=1}^{\infty}
$
where $\lambda_j =j+\mathrm{ln}^2\, j,$ $j=1,2,\dots$,

\noindent
$
2) \ \mathcal M=\{ \mu_j\}_{j\in\mathbb Z'},
$
where
$\mu_j =j+\mathrm{ln}^2\, |j|$,
$j\in\mathbb Z'=\mathbb Z\setminus\{ 0\}$,
satisfy (\ref{PW-1}). 

\noindent
But (\ref{PW})
fails for each of them.
}

\smallskip

{\sl The work is supported by Ministry  of Science and  Higher education of Russian Federation (code of scientific theme FZWU-2020-0027).}

\end{document}